\documentclass{article}
\usepackage{amsmath,amsthm,amsfonts,amssymb}
\usepackage[dvips]{graphicx}
\usepackage[english,french]{babel}

\newtheorem{e-proposition}[theorem]{Proposition}

\newtheorem{e-definition}[theorem]{Definition\rm}

\def\og{\leavevmode\raise.3ex\hbox{$\scriptscriptstyle\langle\!\langle$~}}
\def\fg{\leavevmode\raise.3ex\hbox{~$\!\scriptscriptstyle\,\rangle\!\rangle$}}

\begin{document}

\selectlanguage{english}
\begin{center}
\textbf{Solutions fortes de l'\'{e}quation de l'\'{e}nergie }
\end{center}
\begin{center}
\text{Ch\'erif Amrouche $^{1}$},
\text{Macaire Batchi  $^{1, 2}$},
\text{Jean Batina  $^{2}$}.
\end{center}
\begin{center}
\textit{1 Laboratoire de Math\'ematiques Appliqu\'ees CNRS UMR 5142}
\\
\textit{2 Laboratoire de Thermique Energ\'etique et Proc\'ed\'es}
\\
\textit{Universit\'e de Pau et des Pays de l'Adour}\\
\textit{Avenue de l'Universit\'e 64000 Pau, France}
\end{center}
\begin{abstract}
\selectlanguage{english} In this paper, we give some existence results of stong solutions 
for the energy equation associated to the Navier-Stokes equations with nonhomogeneous boundary conditions 
in two dimension.
\end{abstract}
%
\smallskip
\noindent\textbf{Mots Cl\'es}: Equation de l'\'energie, \'equations
Navier-Stokes, fluide incompressible, convection forc\'ee,
conditions aux limites non homog\`enes, temp\'erature de paroi.


\section{Introduction}

Le probl\`{e}me thermique est d\'{e}crit par l'\'{e}quation de
l'\'{e}nergie avec la convection assur\'{e}e par la vitesse du
fluide du syst\`{e}me (2).\ La temp\'{e}rature $\theta\left(
\boldsymbol{x}\mathbf{,}t\right) $ du fluide satisfait le
syst\`{e}me suivant, dans lequel on suppose qu'il n'y a
pas de source ext\'{e}rieure de chaleur :%
\begin{equation}
\left\{
\begin{array}{lll}
\dfrac{\partial \theta}{\partial t}+(\boldsymbol{v}.\nabla
)\theta-a\triangle \theta=0\text{
\quad } & \text{dans\quad } & Q_{T}=\Omega \times \left] 0,T \right[ , \\
\theta=\theta_{\infty } & \text{sur } & \Gamma _{0}\times \left] 0,T\right[ , \\
\theta=\theta_{p} & \text{sur } & \Gamma _{2}\times \left] 0,T\right[ \\
\dfrac{\partial \theta}{\partial x}=0 & \text{sur } & \Gamma _{1}\times \left] 0,T%
\right[ , \\
\theta\left( 0\right) =\theta_{0} & \text{dans} & \Omega .%
\end{array}%
\right.
\end{equation}%
avec $a>0\medskip $,\quad $T>0$,\quad $\theta_{\infty}$,\quad
$\theta_{p}$ \quad r\'eels donn\'es et $\theta_{0}$ donn\'ee. En
outre, on a le champ de vitesses $\boldsymbol{v}$ solution des
\'equations de Navier-Stokes:
\begin{equation}
\left\{
\begin{array}{lll}
\dfrac{\partial \boldsymbol{v}}{\partial t}-\nu \triangle \boldsymbol{v}%
\text{ }\mathbf{+}\text{ }\boldsymbol{v}\mathbf{.\nabla
}\boldsymbol{v}\text{ }\mathbf{+}\text{ }\mathbf{\nabla }p=0\quad &
\text{dans\quad } &
Q_{T}=\Omega \times \left] 0,T\right[ , \\
\text{div }\boldsymbol{v}\text{ }\mathbf{=}\text{ }0 & \text{dans} &
Q_{T},
\\
\boldsymbol{v}\text{ }\mathbf{=}\text{ }\boldsymbol{g} & \text{sur}
& \Sigma
_{T}=\Gamma \times \left] 0,T\right[ , \\
\boldsymbol{v}\mathbf{(}0\mathbf{)=}\text{ }\boldsymbol{v}_{0} &
\text{dans}
& \Omega .%
\end{array}%
\right.
\end{equation}

\smallskip

\noindent o\`{u} $\boldsymbol{g}$ , $\boldsymbol{v}_{0}$ et $T>0$ sont donn%
\'{e}s. On suppose que :%
\begin{equation}
\begin{array}{ll}
\text{div }\boldsymbol{v}_{0}\text{ }\mathbf{=}\text{ }0\text{ }\
\text{dans
\quad }\Omega , & \boldsymbol{v}_{0}.\boldsymbol{n}=0\text{ \quad sur \quad }%
\Gamma ,%
\end{array}
\end{equation}%
et
\begin{equation}
\boldsymbol{g}.\boldsymbol{n}=0\text{ \quad sur \quad }\Sigma _{T}.
\end{equation}
\\
\noindent On suppose ici que $\Omega$ est un ouvert born\'e de
classe C$^{1,1}$.
\smallskip
\noindent On se ram\`{e}ne aux variables adimensionnelles o\`{u}
l'on pose pour simplifier :

\begin{equation*}
\theta ^{\ast }=\dfrac{\theta-\theta_{\infty%
}}{\theta_{p}-\theta_{\infty }}.
\end{equation*}

\smallskip

\noindent De telle sorte que le syst\`{e}me (1) devient :

\begin{equation}
\left\{
\begin{array}{lll}
\dfrac{\partial \theta ^{\ast }}{\partial t}+(\boldsymbol{v}.\nabla )\theta
^{\ast }-a\triangle \theta ^{\ast }=0\text{ \ \quad } & \text{dans\quad } &
\Omega \times \left] 0,T\right[ , \\
\theta ^{\ast }=0\text{ } & \text{sur } & \Gamma _{0}\times \left] 0,T\right[
, \\
\theta ^{\ast }=1 & \text{sur } & \Gamma _{2}\times \left] 0,T\right[ , \\
\dfrac{\partial \theta ^{\ast }}{\partial x}=0 & \text{sur } & \Gamma
_{1}\times \left] 0,T\right[ , \\
\theta ^{\ast }\left( 0\right) =\theta _{0}^{\ast }\text{ } & \text{dans} &
\Omega .%
\end{array}%
\right.
\end{equation}%
avec $\theta _{0}^{\ast }=\dfrac{\theta_{0}-\theta_{\infty }}{\theta_{p}-\theta_{\infty }}%
.\medskip $

\noindent Pour \'{e}tudier le probl\`{e}me (5), on se ram\`{e}ne
\`{a} des conditions aux limites homog\`{e}nes en posant

\begin{equation*}
\tilde{\theta} =\theta ^{\ast }-\theta _{s}\text{ }
\end{equation*}

\smallskip

\noindent o\`{u} $\theta _{s}$ repr\'{e}sente la fonction de
temp\'{e}rature \`{a} la paroi  et v\'{e}rifie
\begin{equation}
\left\{
\begin{array}{lll}
\triangle \theta _{s}=0\text{\ } & \text{dans} & \Omega , \\
\theta _{s}=0\text{ } & \text{sur } & \Gamma _{0}, \\
\theta _{s}=1\text{ } & \text{sur } & \Gamma _{2}, \\
\dfrac{\partial \theta _{s}}{\partial x}=0\text{ \ } & \text{sur } & \Gamma
_{1}.%
\end{array}%
\right.
\end{equation}

\smallskip

\noindent L'ouvert $\Omega$ \'etant r\'egulier, le probl\`{e}me (6)
poss\`{e}de une solution unique $\theta _{s}\in \underset{1\leq
\text{ }p\text{ }<2}{\cap }W^{1,p}\left( \Omega \right) .$

\noindent Alors (5) s'\'{e}crit%
\begin{equation}
\left\{
\begin{array}{lll}
\dfrac{\partial \theta }{\partial t}+(\boldsymbol{v}.\nabla )\theta
-a\triangle \theta =-(\boldsymbol{v}.\nabla )\theta _{s}\quad \text{ } &
\text{dans\quad } & Q_{T}=\Omega \times \left] 0,T\right[ , \\
\theta =0\text{ \ \ } & \text{sur } & \Gamma _{0}\times \left] 0,T\right[ ,
\\
\theta =0\text{ } & \text{sur } & \Gamma _{2}\times \left] 0,T\right[ , \\
\dfrac{\partial \theta }{\partial x}=0\text{ } & \text{sur } & \Gamma
_{1}\times \left] 0,T\right[ , \\
\theta \left( 0\right) =\theta _{0} & \text{dans} & \Omega .%
\end{array}%
\right.
\end{equation}%
avec%
\begin{equation*}
\theta _{0}=\theta _{0}^{\ast }-\theta _{s},
\end{equation*}%
\noindent et o\`u pour simplifier l'\'ecriture, on a not\'e $\theta$
au lieu de $\tilde{\theta}$ dans (7).
\bigskip

\noindent \textbf{Remarque 1.1\medskip }

\noindent Comme $\theta _{s}\in \underset{1\leq p<2}{\cap }W^{1,p}\left(
\Omega \right) $ et $\ \boldsymbol{v}$ $\mathbf{\in }$ $L^{2}\mathbf{(}0,T%
\mathbf{;\mathbf{H}}^{2}(\Omega ))\smallskip $

\noindent alors\smallskip

$\ \ \ \ \ \ \ \ \ \boldsymbol{v}\mathbf{.\nabla }\theta _{s}$ $\mathbf{\in }
$ $L^{2}\mathbf{(}0,T\mathbf{;}L^{p}\left( \Omega \right) )$ $\forall 1\leq
p<2.\square \medskip \bigskip $

\noindent Afin de r\'{e}soudre le probl\`{e}me (7), nous allons
utiliser \`{a} nouveau la m\'{e}thode de Galerkin. Pour cela, on
d\'efinit l'espace $\Phi =\left\{ \varphi \in H^{1}(\Omega )\text{, }%
\varphi =0\text{ sur }\Gamma _{0}\cup \Gamma _{2}\right\} .\medskip
$ \noindent Le lemme qui suit permet d'obtenir une base sp\'{e}ciale
adapt\'{e}e.\medskip

\noindent \textbf{Lemme 1.2 }\textit{\ Il existe une suite }$(\psi
_{j})_{j\geq 1}$\textit{\ de }$\Phi $\textit{\ et une suite
}$(\lambda _{j})_{j\geq 1}$\textit{\ de r\'{e}els tels que :\medskip
}

$\hspace{1.5cm}\hspace{1.5cm}\lambda _{j}>0,$ \textit{\ }$\underset{%
j\rightarrow \infty }{\lim }\lambda _{j}=+\infty, \medskip $

$\hspace{1.5cm}\hspace{1.5cm}\forall \varphi \in \Phi $\textit{, \
}$\left( \left( \psi _{j},\varphi \right) \right) =\lambda
_{j}\left( \psi _{j}\mid \varphi \right), \medskip $

$\hspace{1.5cm}\hspace{1.5cm}\left( \psi _{j}\mid \psi _{k}\right) =\delta
_{jk},$\textit{\ \ }$\left( \left( \psi _{j},\varphi \right) \right)
=\lambda _{j}\delta _{jk}.\medskip \bigskip $

\noindent \textbf{D\'{e}monstration.} Pour $f\in L^{2}\mathbf{(}(\Omega )%
\mathbf{)}$ , soit $\psi \in H^{2}(\Omega )$ l'unique solution de%
\begin{equation*}
\left\{
\begin{array}{lll}
-\triangle \psi =f\text{ } & \text{dans} & \Omega , \\
\text{\ \ }\psi \text{ }\mathbf{=}\text{ }0\text{ } & \text{sur } & \Gamma
_{0}\cup \text{ }\Gamma _{2}, \\
\dfrac{\partial \psi }{\partial x}=0 & \text{sur } & \Gamma _{1}.%
\end{array}%
\right.
\end{equation*}

\noindent L'op\'{e}rateur $\wedge $ $:\qquad $%
\begin{equation*}
\begin{array}{ll}
f\text{ }\longmapsto \text{ }\psi & \text{lin\'{e}aire et continue}, \\
L^{2}(\Omega )\mathbf{\longrightarrow }H^{2}(\Omega ) &
\end{array}%
\end{equation*}%
\noindent est consid\'{e}r\'{e} comme op\'{e}rateur de $L^{2}(\Omega
)$ dans lui m\^{e}me, $\wedge $ est compact.\ De plus, il est
auto-adjoint :\smallskip

$\qquad
\begin{array}{ll}
\left( \wedge f_{1},f_{2}\right)  =\left( u_{1},-\triangle%
u_{2}\right) =\ \ \displaystyle\int_{\Omega }\nabla u_{1}.\nabla u_{2}dx%
=\ \left( f_{1},\wedge f_{2}\right) .%
\end{array}%
$

$\ \ \ \ \ \ \ \ $\ \ \ \ \ \ \ \ \ \ \ \ \ \ \ \ \ \ \ \ \

\noindent Par cons\'{e}quent, $L^{2}(\Omega )$ poss\`{e}de une base
hilbertienne form\'{e}e de vecteurs propres de $\wedge \medskip $

$\qquad
\begin{array}{lll}
\wedge \psi _{j}=\mu _{j}\psi _{j}, & \mu _{j}\in
\mathbb{R}
, & \ \underset{j\longrightarrow \infty }{\mu _{j}\longrightarrow 0}%
\end{array}%
$

$\ \ \ \ \ \ \ $

\noindent On a donc\medskip

$\ \ \ \ \ \ \ \ \left( \psi _{i}\mid \psi _{j}\right) =\delta _{ij}\medskip
$

\noindent et pour tout $\varphi \in H^{1}(\Omega )$ tel que $\varphi =0$ sur
$\Gamma _{0}$ $\cup $ $\Gamma _{2},$ on a :\medskip

$\qquad
\begin{array}{l}
\left( \wedge \psi _{i}\mid \varphi \right) =\mu _{i}\left( \psi _{i}\mid
\varphi \right)%
\end{array}%
$

\noindent \textit{i.e.}%
\begin{equation*}
\left\{
\begin{array}{lll}
-\triangle \psi _{j}=\dfrac{1}{\mu _{j}}\psi _{j}\text{ } & \text{dans} &
\Omega , \\
\text{\ }\psi _{j}\text{ }\mathbf{=}\text{ }0 & \text{sur } & \Gamma
_{0}\cup \text{ }\Gamma _{2}, \\
\dfrac{\partial \psi _{j}}{\partial x}=0 & \text{sur } & \Gamma _{1}.%
\end{array}%
\right.
\end{equation*}
\noindent Notons que
\begin{equation*}
\left( \psi _{j}\mid \psi _{k}\right) =\lambda _{j}\delta _{jk}\text{ \ o%
\`{u}}\ \text{\ }\lambda _{j}=\dfrac{1}{\mu _{j}}
\end{equation*}


\noindent et que $\left( \psi _{j}\right) _{j\geq 1}$ est une base
orthogonale de l'espace $\Phi.\medskip $

\noindent En particulier pour tout $\varphi \in \Phi \mathbf{,}$ il existe $%
\varphi _{m}\in \left\langle \psi _{1},...,\psi _{m}\right\rangle $ tel que $%
\varphi _{m}\rightarrow \varphi $ dans $H^{1}(\Omega )$ lorsque $%
m\rightarrow \infty .\square \medskip $

\section{Existence de solution}

Consid\'{e}rons maintenant le probl\`{e}me auxiliaire suivant :%
\begin{equation}
\left\{
\begin{array}{lll}
\dfrac{\partial \theta }{\partial t}+(\boldsymbol{v}.\nabla )\theta
-a\triangle \theta =h\text{ } & \text{dans} & \Omega \times \left] 0,T\right[
, \\
\theta =0\text{ \ \ } & \text{sur } & \left( \Gamma _{0}\cup \Gamma
_{2}\right) \times \left] 0,T\right[ , \\
\dfrac{\partial \theta }{\partial x}=0\text{ } & \text{sur } & \Gamma
_{1}\times \left] 0,T\right[ , \\
\theta \left( 0\right) =\theta _{0} & \text{dans} & \Omega .%
\end{array}%
\right.
\end{equation}

\medskip

\noindent o\`{u} $h$ et $\theta _{0}$ donn\'{e}s.\bigskip

\noindent On se propose de d\'{e}montrer le r\'{e}sultat suivant\bigskip
\medskip

\noindent \textbf{Lemme 2.1}\textit{\ \ Si }$\theta _{0}\in
H^{1}\left( \Omega \right) $ \textit{et }$\theta _{0}=0$ \textit{sur
}$\Gamma _{0}\cup \Gamma _{2}$ \textit{et si} $h\in L^{2}\left(
0,T;L^{2}\left( \Omega \right) \right) ,$ \textit{alors le
probl\`{e}me (8) poss\`{e}de une solution unique telle que}
\begin{equation*}
\theta \in L^{2}(0,T;H^{2}\left( \Omega \right) )\cap L^{\infty
}(0,T;H^{1}\left( \Omega \right) )
\end{equation*}%
\begin{equation*}
\theta ^{\prime }\in L^{2}(0,T;L^{2}(\Omega )).\text{ \ \ \ \ \ \ }
\end{equation*}%
\noindent \textbf{D\'{e}monstration.\smallskip} \noindent  On
d\'{e}finit $\theta _{m}\left( t\right) $ une solution
approch\'{e}e du probl\`{e}me (8) par :%
\begin{equation*}
\ \ \ \theta _{m}\left( t\right) =\sum_{j=1}^{m}g_{jm}\left(
t\right) \psi _{j}
\end{equation*}
\noindent De sorte que
\begin{equation}
\left( \theta _{m}^{\prime }\left( t\right) ,\psi _{j}\right)
+\alpha \left( \left( \theta _{m}\left( t\right) ,\psi _{j}\right)
\right) +b\left( \boldsymbol{v},\theta _{m}\left( t\right) ,\psi
_{j}\right) =\left( h,\psi _{j}\right)
\end{equation}


\noindent avec%
\begin{equation}
\left\{
\begin{array}{lll}
\theta _{m}\left( 0\right) =\theta _{0m}\in \left\langle \psi _{1},...,\psi
_{m}\right\rangle ,\\ \theta _{0m}\rightarrow \theta_{0} & \text{dans }%
H^{1}\left( \Omega \right) & \text{lorsque \ }m\rightarrow \infty,\\
\theta _{0}  =0\text{ sur }\Gamma_{0.}%
\end{array}%
\right.
\end{equation}

\noindent \textit{i) Estimation 1\medskip }

Multiplions (9) par $g_{jm}\left( t\right) $ et sommons sur $j$ de
$1$ \`{a} $m$ :

\begin{equation*}
\begin{array}{ll}
\dfrac{1}{2}\dfrac{d}{dt}\left\vert \theta _{m}\left( t\right) \right\vert
^{2}+\alpha \left\Vert \theta _{m}\left( t\right) \right\Vert ^{2} & =\left(
h,\theta _{m}\left( t\right) \right) \\
& \leq C\left\vert h\left( t\right) \right\vert \text{ }\left\Vert \theta
_{m}\left( t\right) \right\Vert%
\end{array}%
\end{equation*}%
$\ \ \ \ \ \ \ \ \ $

\begin{equation}
\dfrac{d}{dt}\left\vert \theta _{m}\left( t\right) \right\vert
^{2}+\alpha \left\Vert \theta _{m}\left( t\right) \right\Vert
^{2}\leq \dfrac{1}{\alpha C^{2}}\left\vert h\left( t\right)
\right\vert ^{2}
\end{equation}%
\smallskip

\noindent En int\'{e}grant (11) de $0$ \`{a} $s,$ o\`{u} $s\in \left] 0,T%
\right] ,$ on obtient\medskip

$\ \ \ \ \ \ \ \ \ \left\vert \theta _{m}\left( t\right) \right\vert
^{2}\leq \left\Vert \theta _{0}\right\Vert ^{2}+\dfrac{1}{\alpha C^{2}}%
\left\Vert h\right\Vert _{L^{2}\left( 0,T;L^{2}\left( \Omega \right) \right)
}^{2},\medskip $

\noindent et on d\'{e}duit \ que\medskip $\ \ \  $%
\begin{equation*}
\ \theta _{m}\in \text{born\'{e} de \ }L^{\infty }(0,T;L^{2}\left( \Omega
\right) ).
\end{equation*}%

\noindent Ensuite, en int\'{e}grant de nouveau (11) entre $0$ et
$T$, il
vient\smallskip $\ \ \ \ \ \ \ \ $%
\begin{equation*}
\ \theta _{m}\in \text{born\'{e} de \ }L^{2}(0,T;H^{1}\left( \Omega
\right) ).
\end{equation*}

\noindent Finalement, on a :

\begin{equation}
\theta _{m}\in \text{born\'{e} de }L^{\infty }(0,T;L^{2}\left(
\Omega \right) )\cap L^{2}(0,T;H^{1}\left( \Omega \right) ).
\end{equation}%

\noindent \textit{ii) Estimation 2\medskip }

Multiplions (9) par $\lambda _{j}g_{jm}\left( t\right) $ puis sommons sur $%
j$ et gr\^{a}ce au choix de la base $\left( \psi _{j}\right) $
d\'{e}finie au lemme 1.2, on a :\smallskip

\bigskip\ $\ \ \ \ \ \ \ \ \ \dfrac{1}{2}\dfrac{d}{dt}\left\Vert \theta
_{m}\left( t\right) \right\Vert ^{2}+\alpha \left\vert \triangle \theta
_{m}\left( t\right) \right\vert ^{2}+b\left( \boldsymbol{v},\theta
_{m}\left( t\right) ,\triangle \theta _{m}\left( t\right) \right) =\left(
h,\triangle \theta _{m}\left( t\right) \right) .$

\noindent D'o\`{u}\smallskip

\bigskip\ $\ \ \ \ \ \ \ \ \dfrac{1}{2}\dfrac{d}{dt}\left\Vert \theta
_{m}\left( t\right) \right\Vert ^{2}+\alpha \left\vert \triangle \theta
_{m}\left( t\right) \right\vert ^{2}\leq \left\vert b\left( \boldsymbol{v}%
,\theta _{m}\left( t\right) ,\triangle \theta _{m}\left( t\right) \right)
\right\vert +\left\vert \left( h,\triangle \theta _{m}\left( t\right)
\right) \right\vert $

\noindent Mais,\medskip

$\qquad
\begin{array}{ll}
\ \left\vert b\left( \boldsymbol{v},\theta _{m}\left( t\right) ,\triangle
\theta _{m}\left( t\right) \right) \right\vert & \leq \left\Vert \boldsymbol{%
v}\right\Vert _{L^{\infty }\left( \Omega \right) }\left\Vert \theta
_{m}\left( t\right) \right\Vert \left\vert \triangle \theta _{m}\left(
t\right) \right\vert \\
& \leq \dfrac{3}{4\alpha }\left\Vert \boldsymbol{v}\right\Vert _{L^{\infty
}\left( \Omega \right) }^{2}\left\Vert \theta _{m}\left( t\right)
\right\Vert ^{2}+\dfrac{\alpha }{4}\left\vert \triangle \theta _{m}\left(
t\right) \right\vert ^{2}.%
\end{array}%
$

$\ \ \ \ \ \ \ \ \ \ \ \ \ \ \ \ \ \ \ \ \ \ \ \ \ \ \ \ \ \ \ \ \ \ \ \ \ \
\ \ \ $

\noindent Et par cons\'{e}quent, on a\medskip

$\qquad
\begin{array}{ll}
\dfrac{1}{2}\dfrac{d}{dt}\left\Vert \theta _{m}\left( t\right) \right\Vert
^{2}+\alpha \left\vert \triangle \theta _{m}\left( t\right) \right\vert
^{2}\leq & \dfrac{3}{4\alpha }\left\Vert \boldsymbol{v}\right\Vert
_{L^{\infty }\left( \Omega \right) }^{2}\left\Vert \theta _{m}\left(
t\right) \right\Vert ^{2}+\dfrac{\alpha }{4}\left\vert \triangle \theta
_{m}\left( t\right) \right\vert ^{2} \\
& +\dfrac{3}{4\alpha }\left\vert h\left( t\right) \right\vert ^{2}+\dfrac{%
\alpha }{4}\left\vert \triangle \theta _{m}\left( t\right) \right\vert ^{2},%
\end{array}%
$

$\smallskip$

\noindent et\medskip

$\ \ \ \ \dfrac{d}{dt}\left\Vert \theta _{m}\left( t\right) \right\Vert
^{2}+\alpha \left\vert \triangle \theta _{m}\left( t\right) \right\vert
^{2}\leq \dfrac{3}{2\alpha }\left\Vert \boldsymbol{v}\right\Vert _{L^{\infty
}\left( \Omega \right) }^{2}\left\Vert \theta _{m}\left( t\right)
\right\Vert ^{2}+\dfrac{3}{2\alpha }\left\vert h\left( t\right) \right\vert
^{2}.\bigskip $

\noindent En utilisant le lemme de Gronwall, puis l'int\'{e}gration entre $0$
et $T,$ on d\'{e}duit que

\begin{equation}
\theta _{m}\in \text{born\'{e} de }L^{2}(0,T;H^{2}\left( \Omega
\right) )\cap L^{\infty }(0,T;H^{1}\left( \Omega \right) ).
\end{equation}%
\noindent \textit{iii) Estimation 3\medskip\ }

Multiplions (9) par $g_{jm}^{\prime }\left( t\right) $ et sommant sur $j$ $%
=$ $1$ \`{a} $m,$ on a :\medskip

$\qquad
\begin{array}{ll}
\left\vert \theta _{m}^{\prime }\left( t\right) \right\vert ^{2}+\dfrac{%
\alpha }{2}\dfrac{d}{dt}\left\Vert \theta _{m}\left( t\right) \right\Vert
^{2} & \leq \left\vert b\left( \boldsymbol{v}\left( t\right) ,\theta
_{m}\left( t\right) ,\theta _{m}^{\prime }\left( t\right) \right)
\right\vert +\left\vert \left( h,\theta _{m}^{\prime }\left( t\right)
\right) \right\vert \\
& \leq \left\Vert \boldsymbol{v}\left( t\right) \right\Vert _{L^{\infty
}\left( \Omega \right) }\left\Vert \theta _{m}\left( t\right) \right\Vert
\left\vert \theta _{m}^{\prime }\left( t\right) \right\vert +\left\vert
h\left( t\right) \right\vert \left\vert \theta _{m}^{\prime }\left( t\right)
\right\vert .%
\end{array}%
$

$\ \ \ \ \ \ \ \ \ \ \ \ \ \ \ \ \ \ \ \ \ \ \ \ \ \ \ \ \ \ \ \ \ \ \ \ $\ $%
\ $\ $\ \ \ \ \ \ \ \ $

\noindent On a donc :\medskip

$\qquad
\begin{array}{ll}
\ \left\vert \theta _{m}^{\prime }\left( t\right) \right\vert ^{2}+\alpha
\dfrac{d}{dt}\left\Vert \theta _{m}\left( t\right) \right\Vert ^{2} & \leq
C\left\Vert \boldsymbol{v}\left( t\right) \right\Vert _{L^{\infty }\left(
\Omega \right) }^{2}\left\Vert \theta _{m}\left( t\right) \right\Vert
^{2}+C\left\vert h\left( t\right) \right\vert ^{2}.%
\end{array}%
$

$\ \ \ \ \ \ \ \ \ $

\noindent En int\'{e}grant entre $0$ et $t,$ on obtient
\begin{eqnarray*}
\displaystyle\int_{0}^{t}\left\vert \theta _{m}^{\prime }\left(
s\right) \right\vert ^{2}ds+\alpha \left\Vert \theta _{m}\left(
t\right) \right\Vert ^{2}\leq C\displaystyle\int_{0}^{t}\left\Vert
\boldsymbol{v}\left( t\right) \right\Vert _{L^{\infty }\left( \Omega
\right) }^{2}\left\Vert \theta _{m}\left( s\right) \right\Vert ^{2}
ds\\
+C\displaystyle\int_{0}^{t}\left\vert h\left( s\right)
\right\vert ^{2} ds+\alpha \left\Vert \theta _{0}\right\Vert
^{2}.\medskip
\end{eqnarray*}
\noindent Il en r\'{e}sulte alors que :
\begin{equation}
\theta _{m}^{\prime }\in \text{born\'{e} de }L^{2}(0,T;L^{2}(\Omega )).
\end{equation}%

\noindent \textit{iv) Passage \`{a} la limite.\medskip\ }

On utilise ici le r\'{e}sultat de compacit\'{e} suivant (voir Temam
$\left[ 7\right] ,$ Lions $\left[ 5\right] $) :\bigskip

\noindent \textbf{Th\'{e}or\`{e}me 2.2} \ \textit{Soit }$B_{0},B,B_{1}$%
\textit{\ trois espaces de Banach avec }
\begin{equation*}
B_{0}\subset B\subset B_{1},B_{0}\text{ \textit{et} }B_{1}\text{ \textit{%
\'{e}tant reflexifs}}
\end{equation*}
\noindent \textit{et on suppose que}
\begin{equation*}
\text{\textit{l'injection de} }B_{0}\rightarrow B\text{ \textit{est compacte}%
.}
\end{equation*}


\noindent \textit{Soit }
\begin{equation*}
W=\left\{ v\in L^{p_{0}}\left( 0,T;B_{0}\right) ,\text{ }v^{\prime }=\dfrac{%
dv}{dt}\in L^{p_{1}}\left( 0,T;B_{1}\right) \right\}
\end{equation*}


\noindent \textit{avec }$1<p_{i}<+\infty ,$\textit{\ }$i=0,1.\medskip $

\noindent \textit{Alors l'injection de }$W$\textit{\ dans }$L^{p_{0}}\left(
0,T;B\right) $\textit{\ est compacte.\bigskip \medskip }

\noindent Nous appliquons maintenant ce th\'{e}or\`{e}me en posant\medskip

$\quad W=\left\{ \theta \in L^{2}\left( 0,T;H^{2}\left( \Omega \right)
\right) ,\text{ }\theta ^{\prime }=\dfrac{d\theta }{dt}\in
L^{2}(0,T;L^{2}(\Omega ))\right\} \medskip $

\noindent L'injection $W\subset L^{2}\left( 0,T;H^{1}\left( \Omega
\right) \right) $ \'{e}tant compacte, donc il existe une sous-suite
de $\theta _{m},$ encore not\'{e}e $\theta _{m}$, telle que lorsque
$m\rightarrow \infty ,$ les convergences suivantes aient lieu :
\begin{eqnarray}
\text{ }\theta _{m} &\rightarrow &\theta ,\text{ dans \quad }%
L^{2}(0,T;H^{2}(\Omega ))\text{ faible}\\ \text{ }\theta _{m}
&\rightarrow &\theta ,\text{ \ dans \quad }L^{\infty
}(0,T;H^{1}(\Omega ))\text{ faible *}\\
\theta _{m} &\rightarrow &\theta ,\text{ \ dans \quad }L^{2}(0,T;H^{1}(%
\Omega ))\text{ fort}\\
\theta _{m}^{\prime } &\rightarrow &\theta ^{\prime },\text{ dans \quad }%
L^{2}(0,T;L^{2}(\Omega ))\text{ faible.}
\end{eqnarray}%

\noindent Passons \`{a} la limite dans (9), on obtient\medskip

$\qquad
\begin{array}{ll}
\forall \psi \in \Phi , & \left( \theta ^{\prime }\left( t\right) ,\psi
\right) +\alpha \left( \left( \theta \left( t\right) ,\psi \right) \right)
+b\left( \boldsymbol{v},\theta \left( t\right) ,\psi \right) =\left( h,\psi
\right)%
\end{array}%
\medskip $

\noindent Par ailleurs, pour tout $t\in \left[ 0,T\right] ,\medskip $

$\qquad
\begin{array}{ll}
\theta _{m}\left( t\right) \rightarrow \theta \left( t\right) & \text{dans }%
H^{1}(\Omega )\text{ faible}%
\end{array}%
\medskip $

\noindent et en particulier\medskip

$\qquad
\begin{array}{ll}
\theta _{0m}=\theta _{m}\left( 0\right) \rightarrow \theta \left( 0\right) &
\text{dans }H^{1}(\Omega )\text{ faible.}%
\end{array}%
\medskip $

\noindent Donc $%
\begin{array}{l}
\theta _{0}=\theta \left( 0\right) .%
\end{array}%
$ Il est clair finalement que gr\^{a}ce \`{a} (15)-(18), on
v\'{e}rifie les conditions aux limites dans (8).\medskip

\section{R\'{e}solution du probl\`{e}me (7)}

De la remarque 1.1, on note que $\boldsymbol{v}\mathbf{.\nabla
}\theta _{s}$ $\mathbf{\in }$
$L^{2}\mathbf{(}0,T\mathbf{;}L^{p}\left( \Omega \right) )$ \ pour
tout $1\leq p<2.$ Il existe $h_{k}\in \mathcal{D}\left( Q_{T}\right)
$ tel que
\begin{equation*}
h_{k}\rightarrow -\boldsymbol{v}\mathbf{.\nabla }\theta _{s}\text{ dans }%
L^{2}\mathbf{(}0,T\mathbf{;}L^{p}\left( \Omega \right) ).
\end{equation*}
\noindent Il existe donc pour chaque $k$ un unique $\theta _{k}\in $ $%
L^{2}(0,T;H^{2}\left( \Omega \right) )\cap L^{\infty
}(0,T;H^{1}\left( \Omega \right) ),$ $\theta _{k}^{\prime }\in $
$L^{2}(0,T;L^{2}\left( \Omega \right) )$ solution de (8) avec second
membre $h_{k}$ et $\theta _{k}\left( 0\right) =$ $\theta _{0}.$

\noindent Par ailleurs,
\begin{equation}
\dfrac{\partial \theta _{k}}{\partial t}-\alpha \triangle \theta _{k}=h_{k}-(%
\boldsymbol{v}.\nabla )\theta _{k}.
\end{equation}

\noindent Il est clair que d'apr\`{e}s (12)

\begin{equation*}
\theta _{k}\in \text{ \ }L^{2}(0,T;H^{2}\left( \Omega \right) )\cap
L^{\infty }(0,T;H^{1}\left( \Omega \right) )
\end{equation*}
\noindent \textit{i.e.}
\begin{equation*}
\boldsymbol{v}\mathbf{.\nabla }\theta _{s}\text{ }\mathbf{\in }\text{ born%
\'{e} de \ }L^{1}\mathbf{(}0,T\mathbf{;}L^{2}\left( \Omega \right)
)\smallskip.
\end{equation*}

\noindent En particulier $\theta _{k}\in $ \ born\'{e} de $%
L^{2}(0,T;H^{1}\left( \Omega \right) )$ et $h_{k}-(\boldsymbol{v}.\nabla
)\theta _{k}\in L^{2}\mathbf{(}0,T\mathbf{;}L^{p}\left( \Omega \right)
).\medskip $

\noindent De sorte qu'en utilisant les r\'{e}sultats sur l'\'{e}quation de
la chaleur et la convergence de $h_{k}$ vers $-\boldsymbol{v}\mathbf{.\nabla
}\theta _{s}$ dans $L^{2}\mathbf{(}0,T\mathbf{;}L^{p}\left( \Omega \right) )$%
, on d\'{e}duit que pour tout \ $1\leq p<2,$
\begin{equation}
\theta_{k} \in \text{born\'{e} de }L^{2}(0,T;W^{2,p}(\Omega )),
\end{equation}
\begin{equation}
\theta_{k}^{\prime } \in \text{born\'{e} de }L^{2}(0,T;L^{p}(\Omega%
)).
\end{equation}

\noindent On peut finalement passer \`{a} la limite pour montrer que le probl%
\`{e}me (7) admet une solution unique $\theta $ telle que
\begin{eqnarray}
\theta \in L^{2}(0,T;W^{2,p}(\Omega )),%
\end{eqnarray}
\begin{eqnarray}
\theta ^{\prime }\in L^{2}(0,T;L^{p}(\Omega )).%
\end{eqnarray}%
\noindent Alors, on peut \'{e}noncer le\medskip

\noindent \textbf{Th\'{e}or\`{e}me 3.1}\quad\textit{Soit }$\theta
_{s}\in \underset{1\leq p<2}{\cap }W^{1,p}\left( \Omega \right) $,
$\theta _{0}\in
H^{1}\left( \Omega \right) $ et $\ \boldsymbol{v%
}$ $\mathbf{\in }$
$L^{2}\mathbf{(}0,T\mathbf{;\mathbf{H}}^{2}(\Omega )).$
\textit{Alors} \textit{le probl\`{e}me (7) admet une solution unique }$%
\theta $ \textit{v\'{e}rifiant, pour tout} $1\leq p<2$
\begin{equation*}
\begin{array}{ll}
\theta \in L^{2}(0,T;W^{2,p}(\Omega )), &\theta ^{\prime }\in
L^{2}(0,T;L^{p}(\Omega )).%
\end{array}%
\end{equation*}
\medskip
\noindent De ce r\'esultat, on d\'eduit le

\noindent\textbf{Corollaire 3.2}\quad\textit{Soit }$\theta_{0}\in
H^{1}\left( \Omega \right) $ et $\ \boldsymbol{v%
}$ $\mathbf{\in }$
$L^{2}\mathbf{(}0,T\mathbf{;\mathbf{H}}^{2}(\Omega )).$
\textit{Alors} \textit{le probl\`{e}me (1) admet une solution unique }$%
\theta $ \textit{v\'{e}rifiant, pour tout} $1\leq p<2$
\smallskip
\begin{equation*}
\begin{array}{ll}
\theta\in L^{2}(0,T;W^{2,p}(\Omega )), & \theta ^{\prime }\in
L^{2}(0,T;L^{p}(\Omega )).%
\end{array}%
\end{equation*}

\smallskip
\noindent \textbf{Remarque 3.3\smallskip }\quad \noindent Si
$\boldsymbol{v}\in $ $L^{\infty }(0,T;\mathbf{H}^{2}\left(
\Omega \right) ),$ on peut montrer que%
\begin{eqnarray*}
\theta \in L^{2}(0,T;W^{2,p}(\Omega )),\qquad \theta ^{\prime } \in
L^{2}(0,T;L^{p}(\Omega )).\square
\end{eqnarray*}

\end{document}